\documentclass{amsart}

\usepackage{amsmath,amsthm,amsfonts,amssymb,bm,graphicx, mathrsfs}

\usepackage
{hyperref}
\hypersetup{colorlinks=true,citecolor=blue,linkcolor=blue,urlcolor=blue,
pdfstartview=FitH }

\theoremstyle{plain}
\newtheorem{theorem}{Theorem}

\newtheorem{lemma}{Lemma}

\newtheorem{cor}[theorem]{Corollary}

\numberwithin{equation}{section}

\theoremstyle{definition}

\renewcommand{\geq}{\geqslant}
\renewcommand{\leq}{\leqslant}

\newcommand{\changed}[1]{{\color{black} #1}}

\usepackage{calc}
\newsavebox\CBox
\newcommand\hcancel[2][0.5pt]{%
  \changed{\ifmmode\sbox\CBox{$#2$}\else\sbox\CBox{#2}\fi%
  \makebox[0pt][l]{\usebox\CBox}%
  \rule[0.5\ht\CBox-#1/2]{\wd\CBox}{#1}}}

\makeatletter
\DeclareRobustCommand\widecheck[1]{{\mathpalette\@widecheck{#1}}}
\def\@widecheck#1#2{%
    \setbox\z@\hbox{\m@th$#1#2$}%
    \setbox\tw@\hbox{\m@th$#1%
       \widehat{%
          \vrule\@width\z@\@height\ht\z@
          \vrule\@height\z@\@width\wd\z@}$}%
    \dp\tw@-\ht\z@
    \@tempdima\ht\z@ \advance\@tempdima2\ht\tw@ \divide\@tempdima\thr@@
    \setbox\tw@\hbox{%
       \raise\@tempdima\hbox{\scalebox{1}[-1]{\lower\@tempdima\box
\tw@}}}%
    {\ooalign{\box\tw@ \cr \box\z@}}}
\makeatother

\begin{document}

\author{Valentin Blomer}
  
\address{Mathematisches Institut, Bunsenstr. 3-5, 37073 G\"ottingen, Germany} \email{vblomer@math.uni-goettingen.de}

 \title{On triple correlations of divisor functions}

\thanks{First author   supported in part  by the  Volkswagen Foundation and  NSF grant 1128155 while enjoying the hospitality of 
the Institute for Advanced Study. The United States Government is authorized to reproduce and distribute reprints notwithstanding any copyright notation herein.}


\begin{abstract} We prove  an asymptotic formula with power saving error term for a certain triple divisor sum.
  \end{abstract}


\setcounter{tocdepth}{2}  \maketitle 

\maketitle

\section{Introduction}

Additive problems with the divisor function $\tau(n)$   have a long history in number theory, and there is  by now a standard machinery to evaluate 
\begin{equation}\label{quadr}
\sum_{n \leq N} \tau(n)\tau(n+h)
\end{equation}
asymptotically   for $N \rightarrow \infty$, with considerable uniformity in $h$. Such sums come up naturally  as off-diagonal terms in the fourth moment of the Riemann zeta-function. 
The strongest results are due to Meurman \cite{Me}, and we refer to \cite{Mo} for a beautiful spectral interpretation and some history of the problem. 

On the other hand, an asymptotic formula for a triple correlation
$$D_h(N) = \sum_{n \leq N} \tau(n) \tau(n-h) \tau(n+h)$$
seems to be an extremely interesting and challenging open problem. Browning   \cite{B} suggests the asymptotic relation
$$D_h(N)  \sim \frac{11}{8} f(h) \prod_p\left(1 - \frac{1}{p}\right)^2 \left( 1 + \frac{2}{p}\right) N(\log N)^3$$
for an explicit multiplicative function $f(h)$ and proves this on average over $h$:
\begin{equation}\label{browning}
\sum_{h \leq H} D_h(N) = \text{main term } + o(HN (\log N)^3),
\end{equation}
provided $H \geq N^{3/4+\varepsilon}$. The method is based on a lattice point argument and essentially elementary. \\

The divisor function can be viewed as a Fourier coefficient of an Eisenstein series, so this problem calls for the application of spectral theory of automorphic forms. As we will see, this is a powerful tool in this situation and makes it possible to reduce the average over $H$ very substantially from $H \geq N^{3/4+\varepsilon}$ to $H \geq N^{1/3+\varepsilon}$. In addition, in contrast to \eqref{browning} the method produces naturally a power saving error term, and it is flexible enough to accommodate for much more general sums. We start with a 
representative example of some arithmetic significance. 

\begin{cor}\label{thm1} Let $W$ be a smooth function with compact support in $[1, 2]$ and Mellin transform $\widehat{W}$. Let $1 \leq H \leq  N/3$ and let $k \geq 2$ be an integer. Then
\begin{equation}\label{asympt}
\begin{split}
&\sum_h W\Bigl(\frac{h}{H}\Bigr) \sum_{N \leq n \leq 2N} \tau_k(n) \tau(n+h)\tau(n-h) \\
&=  \widehat{W}(1) H N Q_{k+1}(\log N) + O\Bigl(N^{\varepsilon}(H^2 +  NH^{1/2} +  N^{3/2} H^{-1/2} +H N^{1-  \frac{1}{k+2}  } )\Bigr),
\end{split}
\end{equation}
where $\tau_k$ denotes the $k$-fold divisor function, $Q_{k+1}$ is a polynomial  (depending only $k$) of degree $k+1$ and leading constant
\begin{equation}\label{euler}
\frac{1}{(k-1)!} \prod_p \left( 1 + \frac{1}{p}\right)\left(1 - \frac{1}{p^{1 + \delta_{p=2}}} \left(1 - \frac{1}{p+1}\right)^k\right),
\end{equation}
and the implied constant in the error term depends on (the Sobolev norms of) $W$, $k$ and $\varepsilon$. 
\end{cor}


The asymptotic formula \eqref{asympt}  is non-trivial (and in fact with a power-saving error term) for 
\begin{equation}\label{range}
  N^{1/3+\varepsilon} \leq H \leq N^{1-\varepsilon},
\end{equation}  
and it is remarkable that this range is independent of $k$. 
   The upper bound, coming from the error term $H^2$, is only for convenience and can easily be removed. Certainly the problem gets harder as $H$ gets smaller. Also the last error term in \eqref{asympt}  can easily be improved.  The bottleneck for the lower bound is the third error term, which is the limit of the automorphic forms machinery; it requires the full force of spectral theory and a fairly delicate analysis of Bessel functions.  \\

A very different approach to correlation sums with divisor functions uses methods of Green and Tao, in particular the nilpotent Hardy-Littlewood method, see \cite{Ma}. While this is capable of treating  correlation sums of arbitrary order, it cannot -- at the present state of knowledge -- produce power saving error terms, and more importantly, it needs at least an average over two long variables, whereas the point of this article is to resemble as much as possible a one-variable situation. Indeed, the sequence $\{n(n+h)(n-h) \mid   N \leq n \leq 2N, h \leq N^{1/3+\varepsilon}\}$ contains $N^{4/3+\varepsilon}$ elements of   size $N^3$, so it is  as dense as a one-variable ``polynomial'' of degree $9/4-\varepsilon$. Therefore even in the special case $k=2$, Corollary \ref{thm1}   goes, for the first time in the literature,  substantially beyond the quadratic case \eqref{quadr}. \\

If desired, the smoothing in the $h$-sum in Theorem \ref{thm1} can be removed as follows: let $0 < \Delta< 1$ be a parameter and let $W_{\Delta}$ be a smooth function that is 1 on $[1, 2]$ and 0 outside of $[1-\Delta, 2+\Delta]$. An inspection of the proof shows that with this $W$ the error term in \eqref{asympt} becomes 
$$O\Bigl(N^{\varepsilon}(H^2 +  \Delta^{-c}(NH^{1/2} +   N^{3/2} H^{-1/2}) +H N^{1-  \frac{1}{k+2}  } )\Bigr)$$
for some absolute constant $c > 0$ (the weakest estimates give $c=8$, but this can be improved easily if necessary), and so we have
\begin{displaymath}
\begin{split}
&\sum_{H \leq h \leq 2H}   \sum_{N \leq   n \leq 2N} \tau_k(n) \tau(n+h)\tau(n-h) \\
& =   H N Q_{k+1}(\log N) + O\Bigl(N^{\varepsilon}(H^2 +  \Delta^{-c}(NH^{1/2} +  N^{3/2} H^{-1/2}) +H N^{1-  \frac{1}{k+2}  }  + \Delta HN )\Bigr),
\end{split}
\end{displaymath}
which upon choosing $$\Delta = \min\left(\frac{1}{2}, \Bigl(\frac{1}{H} + \frac{N}{H^3}\Bigr)^{1/(2(c+1))}\right)$$ is still non-trivial in the range \eqref{range}.\\

Corollary \ref{thm1} features the $k$-fold divisor function,  but a similar result can be obtained with any reasonable arithmetic  function in place of $\tau_k$. What we really prove, is the following asymptotic formula for general \emph{linear forms} in binary additive divisor problems. 

\begin{theorem}\label{thm2}  Let $W$ be a smooth function with compact support in $[1, 2]$ with Mellin transform $\widehat{W}$. Let $1 \leq H \leq  N/3$ and let $k \geq 2$ be an integer. Let $a_n$, $N \leq n \leq 2N$, be any sequence of complex numbers and let $r_d(n)$ denote the Ramanujan sum. Then
\begin{displaymath}
\begin{split}
\sum_h W\left(\frac{h}{H}\right) \sum_{N \leq n \leq 2N} &a(n) \tau(n+h)\tau(n-h) = H \widehat{W}(1)   \sum_{N \leq n\leq 2N} a_n \sum_d \frac{r_d(2n)}{d^{2}} (\log n + 2\gamma - 2 \log d)^2  \\
& \quad\quad + O\left(   N^{\varepsilon} \left(\frac{H^2}{N^{1/2}} + HN^{1/4} + (HN)^{1/2} + \frac{N}{H^{1/2}}\right)\| a \|_2\right),
\end{split}
\end{displaymath}
where the $O$-constant depends on (the Sobolev norms of) $W$ and $\varepsilon$. 
\end{theorem}

We treat the short $h$-sum as a highly unbalanced additive divisor problem. In contrast to the usual version \eqref{quadr}, in the sum
\begin{equation}\label{shifted}
\sum_{h} W\left(\frac{h}{H}\right)\tau(n+h)\tau(n-h)
\end{equation}
 both arguments $n-h$ and $n+h$ are very big (relative to the length of summation) and highly localized. It is therefore not clear a priori to what extent the analysis resembles the classical situation \eqref{quadr}. As will be apparent from the proof, the expression \eqref{shifted} has also a spectral decomposition, similarly as in \cite{Mo}, but here only the holomorphic cusp forms contribute in a substantial fashion, whereas the contribution of Maa{\ss} forms and Eisenstein series is very small. This is somewhat reminiscent of the analysis in \cite{ST}. 
 
 Even though the starting point of the proof of Theorem \ref{thm2} is  to open one of the divisor functions in \eqref{shifted},   using Jutila's circle method and arguing as in \cite{BM}, one can show an analogous  result for Fourier coefficients of cusp forms in place of $\tau(n)$ (without main term, of course). \\

We end this introduction by reminding the reader of the usual $\varepsilon$-convention that will be applied in this paper. In addition, all implied constants may depend on $\varepsilon$, although this will not be displayed explicitly. We write $A \asymp B$ to mean $A \ll B$ and $B \ll A$. 

\section{Analytic technicalties}

In this section we state some known formulas and estimates for future reference. We start with the Voronoi summation formula for the divisor function (see e.g.\ \cite[Theorem 1.7]{Ju1}). 

\begin{lemma}[Voronoi summation]\label{Voronoi}    
  Let $c$ be a positive integer and $a$ an integer coprime to $c$, and let $W$ be a smooth function compactly supported in $(0,\infty)$.  Then
 \begin{equation*}
 \begin{split}
 & \sum_{n\geq 1} \tau(n)W(n)e\Bigl(\frac{an}{c}\Bigr) \\
 &= 
  \frac{1}{c}  \int_0^{\infty}{(\log
    x+2\gamma-2\log c)W(x)dx}
+  \frac{1}{c} \sum_{\pm} \sum_{n\geq 1}\tau(n)
  e\Bigl(\mp\frac{\overline{a}n}{c}\Bigr)\int_0^\infty W(u)\mathcal{J}_{\pm}\Bigl(\frac{4\pi \sqrt{ nu}}{c}\Bigr) d
  u, 
  \end{split}
\end{equation*} 
 where $\gamma$ is Euler's constant and 
$ \mathcal{J}_+(x) =  -2\pi  Y_{0}(x)$, $\mathcal{J}_-(x) = 4   K_{0}(x).$
 \end{lemma}

The Bessel $K$-function is rapidly decaying for large $x$:
\begin{equation}\label{besselK}
K_0(x) \ll (1 + \log |x|) e^{-x}.
\end{equation}

Next we state the Kuznetsov formula in the notation of \cite{BHM}. We define the following integral transforms for a smooth function $\phi : [0, \infty) \rightarrow \Bbb{C}$ satisfying $\phi(0) = \phi'(0) = 0$, $\phi^{(j)}(x) \ll (1+x)^{-3}$ for $0 \leq j \leq 3$:
\begin{displaymath}
\dot{\phi}(k) = 4i^k \int_0^{\infty} \phi(x) J_{k-1}(x) \frac{dx}{x}, \quad \tilde{\phi}(t) = 2\pi i \int_0^{\infty} \phi(x) \frac{J_{2it}(x) - J_{-2it}(x)}{\sinh(\pi t)} \frac{dx}{x}.\\
\end{displaymath}
We let $\mathcal{B}_k$ be an orthonormal basis of the space of holomorphic cusp forms of level 1 and weight $k$, and we write the Fourier expansion of $f\in\mathcal{B}_k$ as
 \begin{displaymath}
  f(z) = \sum_{n \geq 1} \rho_f(n) (4\pi n)^{k/2} e(nz).
\end{displaymath}
Similarly, for Maa{\ss}  forms $f$ of level $1$ and spectral parameter $t$ we write
\begin{equation*}
  f(z) = \sum_{n \not= 0} \rho_f(n) W_{0, it}(4\pi |n|y) e(nx),
\end{equation*}
where $W_{0, it}(y) = (y/\pi)^{1/2} K_{it}(y/2)$ is a Whittaker function. We fix an orthonormal basis $\mathcal{B}$ of Hecke-Maa{\ss} eigenforms. Finally, we write the Fourier expansion of the (unique) Eisenstein series 
$E(z, s)$  of level 1  at $s = 1/2 + it$   as 
\begin{displaymath}
  E (z, 1/2 + it) =  y^{1/2+it} + \varphi(1/2 + it) y^{1/2-it} + \sum_{n \not= 0} \rho(n, t) W_{0, it}(4\pi |n|y) e(nx). 
\end{displaymath}
Then the following spectral sum formula holds.

\begin{lemma}[Kuznetsov formula]\label{kuznetsov}  Let $\phi$ be as in the previous paragraph, and let $a, b > 0$   be integers. Then
\begin{multline*}
  \sum_{  c\geq 1} \frac{1}{c}S(a, b; c) \phi\left(\frac{4\pi \sqrt{ab}}{c}\right) =  \sum_{\substack{k \geq 2\\ k \text{ even}}} \sum_{f \in \mathcal{B}_k } \dot{\phi}(k) \Gamma(k) \sqrt{ab}  {\rho_f(a)} \rho_f(b)\\
  + \sum_{f \in \mathcal{B} } \tilde{\phi}(t_f) \frac{
    \sqrt{ab}}{\cosh(\pi t_f)} {\rho_f(a)} \rho_f(b) + \frac{1}{4\pi }
  \int_{-\infty}^{\infty}\tilde{\phi}(t) \frac{ \sqrt{ab}}{\cosh(\pi
    t)} {\rho(a, t)} \rho (b, t) dt.
\end{multline*}
\end{lemma}

Finally we need the spectral large sieve inequalities of Deshouillers-Iwaniec  \cite[Theorem 2]{DI}. As we will see later, we only need the inequality for the holomorphic forms, but for convenience we treat all three parts of the spectrum equally.

\begin{lemma}[Spectral large sieve]\label{largesieve} Let $T, M \geq
  1$, and let $(a_m)$, $M \leq m \leq 2M$, be a sequence of complex
  numbers. Then all three quantities
\begin{gather*}
  \sum_{\substack{2 \leq k \leq T\\ k \text{ even}}}\Gamma(k) \sum_{f
    \in \mathcal{B}_k }\Bigl| \sum_m a_m \sqrt{m} \rho_f(m)\Bigr|^2,
 \quad\quad \sum_{\substack{f \in \mathcal{B} \\ | t_f| \leq T}
  }\frac{1}{\cosh(\pi t_f)} \Bigl| \sum_m a_m \sqrt{m} \rho_f(\pm
  m)\Bigr|^2,
  \\
  \int_{-T}^T \frac{1}{\cosh(\pi t)} \Bigl| \sum_m a_m \sqrt{m}
  \rho(\pm m, t)\Bigr|^2 dt
\end{gather*}
are bounded by
\begin{displaymath}
  M^{\varepsilon}  (T^2 + M ) \sum_{m} |a_m|^2. 
\end{displaymath}
\end{lemma}

Next we need to analyze a certain Bessel transform.

\begin{lemma}\label{analyzeW} Let $W$ be a  fixed smooth  function with support in $[1/2, 2]$ satisfying $W^{(j)}(x) \ll_j 1$ for all $j$.   
Fix $C \geq 1$.  For $z \geq 4w  > 0$ and  $z \gg 1$  define 
\begin{displaymath}
  W^{\ast}(z, w) =   \int_0^{\infty} W(y) Y_{0}(4\pi\sqrt{z - yw}) dy. 
\end{displaymath}
Then   we have 
\begin{equation*}
  W^{\ast}(z, w) = z^{-1/4} W_+(z, w) e(2 \sqrt{z}) + z^{-1/4} W_-(z, w) e(-2 \sqrt{z})  + {\rm O}_A(C^{-A})
\end{equation*}
for suitable functions $W_{\pm}$  satisfying 
\begin{equation}\label{boundWpm}
   z^i w^j \frac{\partial^i}{\partial z^i}    \frac{\partial^j}{\partial w^j} W_{\pm}(z, w) \begin{cases}
   = 0, & \sqrt{z}/w \leq C^{-\varepsilon},\\
   \ll_{i, j} C^{\varepsilon(i+j)}, & \text{otherwise.}\end{cases}
 \end{equation}
 for any   $i, j \in \mathbb{N}_0$.  Moreover, if $w$ is smooth function with support in a rectangle $ [c_1, c_2] \times [c_1, c_2]$ for two constants $ c_2 > c_1 > 0$ and ${\tt Z} \gg 1$, ${\tt W} > 0$ are two parameters such that $c_1{\tt Z} \geq 4c_2 {\tt W}$, then the double Mellin transform
 $$\widehat{W}_{\pm}(s, t) := \int_0^{\infty} \int_0^{\infty} W_{\pm}(z, w) w\left(\frac{z}{{\tt Z}}, \frac{w}{{\tt W}}\right) z^{s-1} w^{t-1} dz\, dw$$
 is holomorphic on $\Bbb{C}^2$ and rapidly decaying on vertical lines (i.e.\ $\ll_{A} C^{\varepsilon} (1 + |s|)^{-A} (1 + |t|)^{-A}$).
\end{lemma}
 
\textbf{Proof.} This   was proved for $J_{\nu}$ instead of $Y_0$ in \cite[Lemma 9 \& subsequent remark]{BM} and for $Y_0$ in \cite[Section 3]{BFKMM}. For convenience we sketch the argument. The bounds in \eqref{boundWpm} follow from inserting  the asymptotic formula \cite[8.451.2]{GR} of $Y_0$ and integrating by parts sufficiently often. The holomorphicity of $\widehat{W}_{\pm}(s, t)$ follows from the the compact support of $w$, and its growth on vertical lines can be bounded by partial integration and \eqref{boundWpm}. \\

Finally we bound the Bessel transforms occurring in the Kuznetsov formula for a special weight function. This is essentially due to Jutila \cite[Lemma 3 \& Remark 1]{Ju2}, see also \cite[Lemma 16]{BM}. 
\begin{lemma}\label{jutila} Let $Z \gg 1$, $\tau \in \Bbb{R}$,  $\alpha \in [-4/5,  4/5]$ and $w$ a smooth compactly supported function.  For $$\phi(z) = e^{\pm i z \alpha} w\left(\frac{z}{Z}\right) \left(\frac{z}{Z}\right)^{i\tau}$$ we have
\begin{equation}\label{bounds}
 \dot{\phi}(k) \ll_A \frac{1+|\tau|}{Z} \left( 1 + \frac{k}{Z}\right)^{-A}, \quad \tilde{\phi}(t) \ll_A \left(1+ \frac{|t| + Z}{1 + |\tau| }\right)^{-A}
\end{equation}
for $t \in \Bbb{R}$, $k \in \Bbb{N}$ and any $A \geq 0$. 
\end{lemma}

\textbf{Remark:} In our application later, $Z$ will be typically relatively large while $\tau$ is quite  small, so that  $\tilde{\phi}$ is essentially negligible, while $\dot{\phi}$ is only negligible once $k$ becomes bigger than $Z$. This explains our earlier remark that   only the holomorphic contribution is relevant in our situation. Nevertheless, for a simple uniform treatment we weaken these bounds and combine them to 
\begin{equation}\label{simple}
\dot{\phi}(k)  \ll_A \frac{1+|\tau|}{Z} \left( 1 + \frac{k}{|\tau| + Z}\right)^{-A}, \quad \tilde{\phi}(t)  \ll _A \frac{1+|\tau|}{Z} \left( 1 + \frac{|t|}{|\tau| + Z}\right)^{-A}.
\end{equation}

\textbf{Proof.} We start with the discussion of  the transform $\dot{\phi}$. If $k \geq c Z$ with $c >0$ sufficiently large (depending on the support of $w$), we use the 
  bound $J_{k}(x) \ll e^{-k/5}$ for  $x \leq k/2$ (cf.\ \cite[8.452.1]{GR} or \cite[(4.8)]{Ra}), so that $\dot{\phi}(k) \ll e^{k/5}$, and the first bound in \eqref{bounds} follows. Let us now assume $k \ll Z$. 
Then we insert the Fourier integral \cite[8.411.1]{GR} getting
$$\dot{\phi}(k) = 4i^k \int_0^{\infty} e^{\pm i x \alpha} w\left(\frac{x}{Z}\right) \left(\frac{x}{Z}\right)^{i\tau}  \frac{1}{\pi}  \int_{-\pi}^{\pi} \cos((k-1)\xi - x \sin \xi) d\xi \frac{dx}{x}.$$
Integrating by parts sufficiently often, we see that the $x$-integral is $$\ll_A \left( 1 + \frac{Z |\sin \xi \pm \alpha| }{(1+|\tau|)} \right)^{-A}.$$
We estimate the $\xi$-integral  trivially and obtain again the first bound in \eqref{bounds}. 



For the estimation of $\tilde{\phi}$ we 
first insert the uniform asymptotic expansion \cite[7.13(17)]{EMOT} 
$$ \frac{J_{2it}(x) - J_{-2it}(x)}{\sinh(\pi t)}  = \Re \left(\exp\Bigl(i (\sqrt{(2t)^2 + x^2} - 2 t \, \text{arcsinh}(2t/x)\Bigr) f_A(t, x)\right) + O_A((x+|t|)^{-A})$$
for any fixed $A > 0$, where
$$x^j \frac{\partial^j}{\partial x^j} f_A(t, x) \ll _{A, j} (|t|+x)^{-1/2} \ll 1.$$
The original error term in \cite{EMOT} is only $O(x^{-A})$, but the stronger error term $O((x+|t|)^{-A})$ (possibly with a different $A$) follows from the power series expansion \cite[8.402]{GR} of $J_{it}(x)$ for $x < |t|^{1/3}$, say. 
We then obtain an oscillatory integral with phase
$$H(x) = \pm \alpha x \pm \big(\sqrt{(2t)^2 + x^2} - 2 t \, \text{arcsinh}(2t/x)\big)$$
with derivatives
$$H'(x) = \pm \alpha \pm \frac{\sqrt{(2t)^2 + x^2}}{x}\gg 1 + \frac{|t|}{x}, \quad H^{(j)}(x) \asymp \frac{t^2}{x^j (|t|+x)} \quad (j\geq 2), $$
since $|\alpha| \leq 4/5$.  The following general integration-by-parts lemma  (\cite[Lemma 8.1]{BKY}) with
$$X = \frac{1}{Z}, \quad Q = Z \asymp b-a, \quad U = \frac{Z}{1+|\tau|}, \quad R = 1 + \frac{|t|}{Z}, \quad Y = 1 + \frac{t^2}{|t|+Z}$$
gives the desired second bound in \eqref{bounds}.

 \begin{lemma} \label{integrationbyparts}
 Let $Y \geq 1$, $X, Q, U, R > 0$, 
and suppose that $w$ 
 is a smooth function with support on some interval $[a, b]$, satisfying
\begin{equation*} 
w^{(j)}(t) \ll_j X U^{-j}.
\end{equation*}
Suppose $H$ 
  is a smooth function on $[a, b]$ such that
\begin{equation}\label{condi1}
 |H'(t)| \gg R, \quad 
H^{(j)}(t) \ll_j Y Q^{-j} \,\, \text{for } j=2, 3, \dots.
\end{equation}
Then 
\begin{equation}\label{intbound}
I = \int_{\Bbb{R}} w(t) e^{i H(t)} dt
 \ll_A (b-a) X \left[(QR/\sqrt{Y})^{-A} + (RU)^{-A}\right]
\end{equation}
for any $A \geq 0$. 
\end{lemma}


\section{Proof of Theorem \ref{thm2}}

Most of the time we will consider the sum
$$S(n) := \sum_h W\left(\frac{h}{H}\right) \tau(n+h)\tau(n-h)$$
for $N \leq n \leq 2N$, where $N \geq 3 H$, and we will postpone the sum over $n$ to the very last moment. Clearly we have 
\begin{displaymath}
\begin{split}
S(n) =  \sum_{h} W(h/H) \sum_{ab = n+h} \tau(2n-ab).
\end{split}
\end{displaymath}
Let $V$ be a function that is constantly 1 on $[0, 1/2]$, 0 on $[2, \infty)$ and satisfies $V(x) + V(1/x) = 2$. Then
\begin{displaymath}
S(n) =  \sum_{h} W(h/H) \sum_{ab = n+h} \tau(2n-ab) V\left(\frac{a}{\sqrt{n+h}}\right).  
\end{displaymath}
By the mean value theorem and the trivial bound $\tau(n) \ll n^{\varepsilon}$, we can write
\begin{equation*}\label{err1}
S(n) =  \sum_{h} W(h/H) \sum_{ab = n+h} \tau(2n-ab) V\left(\frac{a}{\sqrt{n}}\right)  + O\left(\frac{H^2}{N^{1-\varepsilon}}\right).
\end{equation*}
Call the main term $S'(n)$. Then
\begin{displaymath}
\begin{split}
S'(n) &= \sum_a V(a/\sqrt{n}) \sum_{m \equiv 2n \, (\text{mod } a)} \tau(m) W\left(\frac{n-m}{H}\right)  \\
& = \sum_a \frac{V(a/\sqrt{n})}{a} \sum_{d \mid a} \underset{b \, (\text{mod }d)}{\left.\sum \right.^{\ast} }  e\left(\frac{-2nb}{d}\right)\sum_m \tau(m) e\left(\frac{mb}{d}\right) W\left(\frac{n-m}{H}\right), 
\end{split}
\end{displaymath}
where the star indicates summation over primitive residue classes. 
We use the Voronoi formula on the $m$-sum. This produces a main term
\begin{equation*}\label{err2}
\begin{split}
& \sum_a \frac{V(a/\sqrt{n})}{a} \sum_{d \mid a}   \frac{r_d(2n)}{d}  \int (\log x + 2\gamma - 2 \log d)W\left(\frac{n-x}{H}\right)   dx\\
& = H \sum_{d, f} V\left(\frac{df}{\sqrt{n}}\right)    \frac{r_d(2n)}{d^2f}  \int (\log (n - Hy) + 2\gamma - 2 \log d)W(y) dy\\
& =  H \sum_{d, f} V\left(\frac{df}{\sqrt{n}}\right)    \frac{r_d(2n)}{d^2f} (\log n + 2\gamma - 2 \log d) \widehat{W}(1) + O\left(\frac{H^2 \log N }{N}\right),
\end{split}
\end{equation*}
where we again applied the mean value theorem in the last step, and 
 $r_d(n)$ denotes the Ramanujan sum. By Mellin inversion, the main term equals
\begin{displaymath}
\begin{split}
&H \widehat{W}(1) \int_{(2)}  \widehat{V}(s) n^{s/2} \zeta(s+1) \sum_d \frac{r_d(2n)}{d^{2+s}} (\log n + 2\gamma - 2 \log d)  \frac{ds}{2\pi i}.\\
\end{split}
\end{displaymath}
We have
$$\sum_{d} \frac{r_d(2n)}{d^{2+s}} = \frac{\sigma_{-1-s}(2n)}{\zeta(s+2)},$$
where $\sigma_s(n) = \sum_{d \mid n} d^s$, and we define more generally
\begin{equation}\label{sigma}
\sigma^{(j)}_{-1}(n) =  (-1)^{j} \frac{d^j}{ds^j} \sigma_{-1-s}(n)|_{s = 0} = \sum_{d \mid n} \frac{(\log d)^j}{d}.
\end{equation}
Shifting the contour to $\Re s = -1+\varepsilon$, this equals
\begin{equation}\label{err3}
\begin{split}
& H \widehat{W}(1)   \sum_d \frac{r_d(2n)}{d^{2}} (\log n + 2\gamma - 2 \log d)^2   + O(H N^{-1/2+\varepsilon})\\
 = & H \widehat{W}(1)   \sum_{i+j \leq 2} c_{i, j} (\log n)^i \sigma_{-1}^{(j)}(2n)   + O(H N^{-1/2+\varepsilon})
\end{split}
\end{equation}
for certain constants $c_{i, j}$ with  
\begin{equation}\label{zeta}
  c_{2, 0} = 1/\zeta(2).
\end{equation}  
   So far we have proved
\begin{equation}\label{sofar}
S(n) =  H \widehat{W}(1)   \sum_d \frac{r_d(2n)}{d^{2}} (\log n + 2\gamma - 2 \log d)^2 + E(n) + O(H^2 N^{\varepsilon-1} + HN^{-1/2 +\varepsilon}),
\end{equation}
where  the error term in the Voronoi formula is given by 
\begin{displaymath}
\begin{split}
E(n) = & \sum_{\pm} \sum_a \frac{V(a/\sqrt{n})}{a} \sum_{d \mid a}  \frac{1}{d} \sum_m \tau(m)  S(-2n, \mp m, d) \int_0^{\infty} \mathcal{J}^{\pm}\left(\frac{4\pi \sqrt{mx}}{d}\right) W\left(\frac{n-x}{H}\right)   dx
\\
 = & H \sum_{\pm} \sum_m \tau(m) \sum_f \frac{1}{f} \sum_d \frac{V(df/\sqrt{n})}{d^2}S(-2n, \mp m, d) \Phi^{\pm}_n(m, d), 
\end{split}
\end{displaymath}
where 
$$\Phi^{\pm}_n(m, d) = \int_0^{\infty}  \mathcal{J}^{\pm}\left(\frac{4\pi \sqrt{m(n - Hy)}}{d}\right) W(y)   dy.$$
By \eqref{besselK} and Weil's bound for Kloosterman sums we can estimate the contribution of the  minus-case trivially by 
\begin{equation}\label{err4}
 \ll H \sum_{m } \sum_{fd \leq 2\sqrt{n}} \frac{\tau(m)(m, d)^{1/2} }{fd^{3/2}} e^{-\sqrt{mn}/d} \ll HN^{-1/4 +\varepsilon}. 
\end{equation}

We proceed to treat the plus-case. This is the most technical part of the argument. First we insert a smooth partition of unity and attach a factor $v(d/D)$ for $D \ll \sqrt{N}/f$, where $v$ is a smooth function with compact support in $[1, 2]$.  Using Lemma \ref{analyzeW} with $C = N$, say, we can write, up to a negligible error, 
$$\Phi^{+}_{n}(m, d) = \sum_{\pm}  \frac{d^{1/2}}{(mn)^{1/4}} \Psi_{\pm}\left( \frac{4\pi \sqrt{mn}}{d}, \frac{4\pi \sqrt{mH}}{d}\right) e\left(\pm 2 \frac{\sqrt{mn}}{d}\right) $$
for certain non-oscillating functions $\Psi_{\pm}$, and we can truncate the $m$-sum, again at the cost of a negligible error, at
\begin{equation}\label{M}
 m \leq M := D^2 N^{1+\varepsilon} H^{-2}.
 \end{equation}
Thus we see that that $E(n)$ is a sum of $O((\log N)^2)$ terms of the form
$$E(n, D, M_0) =  \sum_f \frac{1}{f}\frac{H}{n^{1/4}D^{1/2}} \sum_{M_0 \leq m < 2M_0}  \frac{\tau(m)}{m^{1/4}} \sum_{d} \frac{S(2n, m, d) }{d}  \Theta_{m, n}\left(\frac{4\pi \sqrt{2mn}}{d}\right)$$
for $M_0 \leq M$ plus an error term of size $O(HN^{-1/4+\varepsilon})$ from \eqref{err4}, where
$$\Theta_{m, n}(z) = \sum_{\pm} e^{\pm i z/\sqrt{2}} \tilde{v}\left(\frac{4\pi \sqrt{2mn}}{zD}\right) V\left(\frac{4\pi \sqrt{2m} f}{z}\right) \Psi_{\pm}\left(\frac{z}{\sqrt{2}}, z \frac{H^{1/2}}{(2n)^{1/2}}\right) $$
with $\tilde{v}(x) = v(x)/x$. 
Notice that the support of $\tilde{v}$ localizes  $$z \asymp Z := \sqrt{M_0N}/D, $$
and we have by the size conditions on $D$ and $M_0$ the bounds $1 \ll Z \leq N^{1+\varepsilon}/H.$  
We remember this by attaching   a redundant weight function $w(z/Z)$ to $\Theta_{m, n}(z)$, where $w$  is 1 on $[1/100, 100]$ and 0 outside $[1/200, 200]$. With the aim of applying Lemma \ref{analyzeW}, we also attach the same redundant weight function a second time to $\Psi_{\pm}$ and define
$$ \tilde{\Psi}_{\pm}\left(\frac{z}{\sqrt{2}}, z \frac{H^{1/2}}{(2n)^{1/2}}\right) =  \Psi_{\pm}\left(\frac{z}{\sqrt{2}}, z \frac{H^{1/2}}{(2n)^{1/2}}\right) w\left(\frac{z}{Z}\right).$$
Now we separate variables by brute force Mellin inversion. This can be done at almost no cost, since all weight functions are non-oscillating. We have
\begin{displaymath}
\begin{split}
 & \tilde{v}\left(\frac{4\pi \sqrt{2mn}}{zD}\right) V\left(\frac{4\pi \sqrt{2m} f}{z}\right)\tilde{ \Psi}_{\pm}\left(\frac{z}{\sqrt{2}}, z \frac{H^{1/2}}{(2n)^{1/2}}\right)\\
  &  = \int_{\mathcal{C}} \widehat{\tilde{v}}(s_1) \widehat{V}(s_2) \widehat{\tilde{\Psi}}_{\pm}(s_3, s_4) \left(\frac{4\pi \sqrt{2mn}}{zD}\right)^{-s_1} \left(\frac{4\pi \sqrt{2m} f}{z}\right)^{-s_2} \left(\frac{z}{\sqrt{2}}\right)^{-s_3} \left(z \frac{H^{1/2}}{(2n)^{1/2}}\right)^{-s_4} \frac{ds}{(2\pi i)^4},
\end{split}
\end{displaymath}
where the fourfold contour $\mathcal{C}$ is over the vertical lines $\Re s_1 =  \Re s_3 = \Re s_4 = 0$, $\Re s_2 = \varepsilon$ (since $\widehat{V}$ has a pole at $s = 0$).  By the rapid decay of $\widehat{\tilde{v}}(s_1) \widehat{V}(s_2) \widehat{\tilde{\Psi}}(s_3, s_4) $ along vertical lines, we can truncate the contours at $\Im s_j \ll N^{\varepsilon}$, at the cost of a negligible error. We call this truncated contour $\tilde{\mathcal{C}}$. This gives
\begin{displaymath}
\begin{split}
  E(n, D, M_0) &= \sum_{\pm}  \int_{\tilde{\mathcal{C}}} \widehat{\tilde{v}}(s_1) \widehat{V}(s_2) \widehat{\tilde{\Psi}}_{\pm} (s_3, s_4) (4\pi)^{-s_1 -s_2} (2Z)^{\frac{1}{2}(-s_1-s_2 + s_3 + s_4)}  \sum_f   \frac{ H^{1 - \frac{s_4}{2}}}{f^{1 + s_2}n^{\frac{1}{4} + \frac{s_1}{2} - \frac{s_4}{2} }D^{\frac{1}{2} + s_1}}\\
  &\times  \sum_{M_0 \leq m < 2M_0}  \frac{\tau(m)}{m^{\frac{1}{4} + \frac{s_1}{2} + \frac{s_2}{2}}} \sum_{d} \frac{S(2n, m, d) }{d}  \Xi_{\pm}\left(\frac{4\pi \sqrt{2mn}}{d} \right)\frac{ds}{(2\pi i)^4} + O(N^{-10}), 
\end{split}
\end{displaymath}
where
$$\Xi_{\pm}(z) =  e^{\pm i z/\sqrt{2}} w\left(\frac{z}{Z}\right) \left(\frac{z}{Z}\right)^{s_1 + s_2 - s_3 - s_4}$$
with $\Im(s_1 + s_2 - s_3 - s_4) \ll N^{\varepsilon}$. 

The $d$-sum is in good shape to apply the Kuznetsov formula. Notice that the Selberg eigenvalue conjecture is known in level 1, so   there are no exceptional spectral parameters in the Maa{\ss} spectrum. The corresponding weight function was analyzed in Lemma \ref{jutila}, which we apply with $\alpha = 1/\sqrt{2}$ and $\tau = \Im(s_1 + s_2 - s_3 - s_4) \ll N^{\varepsilon}$. The holomorphic part is given by
$$\sum_{\substack{k \geq 2\\ k \text{ even}}} \sum_{f \in \mathcal{B}_k } \dot{\Xi}_{\pm}(k) \Gamma(k) \sqrt{2nm}  {\rho_f(2n)} \rho_f(m),$$
and similar expressions hold for the Eisenstein and Maa{\ss} spectrum.  It is now an appropriate point to include the sum over $n$. We use the Cauchy-Schwarz inequality and estimate the $s_1, \ldots, s_4$-integrals trivially. This gives
\begin{displaymath}
\begin{split}
  \sum_{N \leq n \leq 2 N} &a(n) E(n, D, M_0)\\
  &\ll N^{\varepsilon} \sup_{\tau_1, \tau_2 \ll N^{\varepsilon}} \sum_f   \frac{ H }{f D^{1/2} } \Biggl(\sum_{\substack{k \geq 2\\ k \text{ even}}} \sum_{f \in \mathcal{B}_k } |\dot{\Xi}_{\pm}(k) |  \Gamma(k) \cdot \Bigl|\sum_{N \leq n \leq 2N} \frac{a_n}{n^{\frac{1}{4} + i\tau_1}} \sqrt{2n} \rho_f(2n)\Bigr|^2\Biggr)^{1/2} \\
  & \quad\quad\times \Biggl(\sum_{\substack{k \geq 2\\ k \text{ even}}} \sum_{f \in \mathcal{B}_k } |\dot{\Xi}_{\pm}(k)|  \Gamma(k)\cdot  \Bigl|\sum_{M_0 \leq m < 2M_0} \frac{\tau_m}{n^{\frac{1}{4} + \frac{\varepsilon}{2}+ i\tau_2}} \sqrt{m} \rho_f(m)\Bigr|^2\Biggr)^{1/2} 
  \end{split}
\end{displaymath}
plus similar Maa{\ss} and Eisenstein terms. 
We apply the spectral large sieve inequalities and \eqref{simple} and  notice that 
 we can truncate the $k$-sum, up to a negligible error, at $k \ll ZN^{\varepsilon}$. This shows
\begin{displaymath}
\begin{split}
 &\sum_{N \leq n \leq 2 N} a(n) E(n, D, M_0)\\
 & \ll N^{\varepsilon}\sum_{f} \frac{H}{f D^{1/2}Z} (Z^2 + N)^{1/2} \Bigl(\sum_{N \leq n \leq 2N} \frac{|a(n)|^2}{n^{1/2}} \Bigr)^{1/2} ( Z^2 + M_0)^{1/2} \Bigl(\sum_{m \asymp M_0} \frac{|\tau(m)|^2}{m^{1/2}} \Bigr)^{1/2}\\
 &\ll   N^{\varepsilon}\sum_{f} \frac{H   D^{1/2}}{f  M_0^{1/4} N^{3/4}} \left(\frac{M_0N}{D^2} + N\right)^{1/2}  \left(\frac{M_0N}{D^2} + M_0\right)^{1/2}   \| a \|_2.
 \end{split}
\end{displaymath}
This is increasing in $M_0$. Inserting the upper bound \eqref{M}, this is at most
$$\ll N^{\varepsilon}\sum_{f} \frac{H^{1/2}   D }{f  N^{1/2}} \left(\frac{ N^2}{H^2} + N\right)^{1/2}  \left(\frac{ N}{D^2} + 1\right)^{1/2}   \| a \|_2.$$
This is non-decreasing in $D$. Inserting  the upper bound $D\ll N^{1/2}/f$, we obtain the bound
$$\ll N^{\varepsilon}\sum_{f \ll N^{1/2} } \frac{H^{1/2}    }{f  } \left(\frac{ N^2}{H^2} + N\right)^{1/2}   \| a \|_2 \ll N^{\varepsilon} \left(\frac{N}{H^{1/2}} + (NH)^{1/2}\right) \| a \|_2.$$
Collecting also the error terms in \eqref{sofar} and \eqref{err4} and summing them over $n$, we obtain Theorem \ref{thm2}. 

\section{Proof of Corollary \ref{thm1}}

By \eqref{err3}, it remains to show that the main term
\begin{equation}\label{remains}
H \widehat{W}(1) \sum_{N \leq n \leq 2N}\sum_{i+j \leq 2} c_{i, j} \tau_k(n) (\log n)^i \sigma_{-1}^{(j)}(2n)
\end{equation}
is, up to an admissible error term, of the desired shape. To this end we study (initially in $\Re  s > 1$) the Dirichlet series
\begin{equation}\label{diri}
L_{i, j}(s) := \sum_n \frac{\tau_k(n) (\log n)^i \sigma_{-1}^{(j)}(2n)}{n^s} = (-1)^i \frac{d^i}{ds^i} L_{0, j}(s).
\end{equation}
By \eqref{sigma} we have $$L_{0, j}(s) = \sum_{d} \frac{(\log d)^j (d, 2)^s}{d^{s+1}} T_{d/(d, 2)}(s),$$ 
where
$$T_{\alpha}(s) = \sum_m \frac{\tau_k(\alpha m)}{m^s} = H_{\alpha}(s) \zeta(s)^k , \quad H_{\alpha}(s) = \sum_{m \mid \alpha^{\infty}} \frac{\tau_k(\alpha m)}{m^s}  \prod_{p \mid \alpha}\left(1 - \frac{1}{p^s}\right)^k,$$
so that
$$L_{i, j}(s) = (-1)^i \frac{d^i}{ds^i} (Z_j(s) \zeta(s)^k), \quad Z_j(s) =  \sum_{d} \frac{(\log d)^j (d, 2)^s}{d^{s+1}} H_{d/(d, 2)}(s).$$
Since $H_{\alpha}(s) \ll_{k,  \sigma_0} \alpha^{\varepsilon}$ in $\Re s \geq \sigma_0 > 0$, the series 
$Z_j(s)$ is holomorphic and uniformly bounded in $\Re s \geq \sigma_0 > 0$. By a standard application of Perron's formula (see e.g.\ \cite[p.\ 133]{Te}) we conclude that
$$\sum_{N\leq n \leq 2N} \tau_k(n) (\log n)^i \sigma_{-1}^{(j)}(2n) = \int_{1+\varepsilon - iT}^{1+\varepsilon+iT} L_{i, j}(s) \frac{(2N)^s - N^s}{s} \frac{ds}{2\pi i} + O\left(\frac{N^{1+\varepsilon}}{T}\right).$$
for $1 \leq T \leq N$. We replace the line of integration by the contour $[1+\varepsilon - iT, \varepsilon - iT] \cup [\varepsilon - iT, \varepsilon + iT]\cup[\varepsilon + iT, 1+\varepsilon +iT]$. Using the convexity bound $\zeta^{(i)}(\sigma + it) \ll_i |t|^{(1-\sigma)/2 + \varepsilon}$ for $\varepsilon \leq \sigma \leq 1+\varepsilon$, $|t| \geq 1$, and picking up to pole at $s=1$, we see that the previous display is 
$$ N P_{k-1+i}(\log N) + O\left(N^{\varepsilon} \Bigl(\frac{N}{T} + T^{k/2}\Bigr)\right) = N P_{k-1+i}(\log N) + O\left(N^{1 - \frac{2}{k+2}+\varepsilon}\right)$$
for a polynomial $P_{k-1+i}$ of degree $k-1+i$ upon choosing $T = N^{2/(k+2)}$.  
(Better error terms could be obtained by arguing more carefully here.)  This shows the formula \eqref{asympt}. 

To compute \eqref{euler}, we recall \eqref{zeta} and observe that only the term $i=2$, $j=0$ in \eqref{remains} is responsible for the leading coefficient. By \eqref{diri} this requires the computation of the Euler product 
$$L_{0, 0}(s)\zeta(s)^{-k}|_{s=1} =  \prod_p \Biggl(\sum_{\alpha=0}^{\infty}  \left(\begin{matrix}\alpha+k-1\\ k-1\end{matrix}\right) \frac{\left(1 - p^{-\alpha-1 - \delta_{p=2}}\right)}{1 - p^{-1}} \frac{1}{p^{\alpha}}\Biggr) \left(1 - \frac{1}{p}\right)^{k}.$$
Using the binomial theorem to evaluate the $\alpha$-sum and incorporating the factor $1/\zeta(2)$ from \eqref{zeta}, we arrive after a short calculation at \eqref{euler}. 
This completes the proof.

\end{document}